\documentclass[11pt]{article}
\usepackage[margin=1in]{geometry}

\usepackage{setspace}
\doublespace

\usepackage{amsmath}
\usepackage{amssymb}
\usepackage[square]{natbib}
\setcitestyle{numbers}

\newtheorem{theorem}{Theorem}[section]

\newtheorem{lemma}{Lemma}[section]
\newtheorem{remark}{Remark}[section]

\newcommand{\Proof}{\textbf{Proof. }}            % Proof
\newcommand{\qed}{$\square$}                     % end of proof

\newcommand{\RR}{\mathbb{R}}
\newcommand{\NN}{\mathbb{N}}
\newcommand{\ZZ}{\mathbb{Z}}
\newcommand{\Zpl}{\ZZ_+}

\newcommand{\ii}{{\mathrm i}}
\newcommand{\ee}{{\mathrm e}}

\newcommand{\dd}{{\mathrm{d}}}

\newcommand{\BI}{\mathrm{Bi}}    % Binomial
\newcommand{\NB}{\mathrm{NB}}   %NBinomial
\newcommand{\G}{\mathrm{G}}    % SCP 2

\newcommand{\Ha}{\mathrm{H}}    % main
\newcommand{\TP}{\mathrm{TP}}   % TP

\newcommand{\norm}[1]{\|#1\|}                    % norm
\newcommand{\ab}[1]{\vert#1\vert}                % absolut value
\newcommand{\Ab}[1]{\Big\vert#1\Big\vert}        % large absol. value

\newcommand{\exponent}[1]{\exp\{#1\}}            % exponent
\newcommand{\Exponent}[1]{\exp\Bigl\{#1\Bigr\}}  % large exponent

\newcommand{\eit}{\ee^{\ii t}}                   % e^{it}

\newcommand{\Expect}{\mathrm{E}}                 % expectation
\newcommand{\Var}{\mathrm{Var}}                  % variance
\newcommand{\qubar}{\overline{q}}
\newcommand{\pbar}{\overline{p}}

\newcommand{\vfi}{\varphi}
\newcommand{\wE}{\w\Expect}
\newcommand{\wexp}{\w\Expect^{+}}
\newcommand{\eL}{{\cal L}}
\newcommand{\w}{\widehat}

\newcommand{\floor}[1]{{\lfloor #1\rfloor}}

\begin{document}
%%%%%%%%%%%%%%%%%%%%%%%%%%%%%%%%%%%%%%%%%%%%%%%%%%%%%%%%%%%%%%%%%%%%%%
%%%%%%%%%%%%%%%%%%%%%%%%%%%%%%%%%%%%%%%%%%%%%%%%%%%%%%%%%%%%%%%%%%%%%%
\title{Non-uniform approximations for sums of discrete m-dependent
random variables}

\author{P. Vellaisamy  and V. \v Cekanavi\v cius    \\
{\small
 Department of Mathematics, Indian Institute of Technology Bombay,} \\
 {\small Powai, Mumbai-
400076, India.}\\{\small  E-mail: pv@math.iitb.ac.in} \\{\small
and}
\\{\small
Department of Mathematics and Informatics, Vilnius University,}\\
{\small Naugarduko 24, Vilnius 03225, Lithuania.}\\{\small E-mail:
vydas.cekanavicius@mif.vu.lt } }

\date{}
%\jobname, \today }

\maketitle
%%%%%%%%%%%%%%%%%%%%%%%%%%%%%%%%%%%%%%%%%%%%%%%%%%%%%%%%%%%%%%%%%%%%%%

\begin{abstract}
Non-uniform estimates are obtained for Poisson, compound Poisson,
translated Poisson, negative binomial and binomial approximations
to sums of of m-dependent integer-valued random variables.
Estimates for Wasserstein metric also  follow easily from our results.  The results are then
exemplified by the approximation of  Poisson binomial distribution, 2-runs and
$m$-dependent $(k_1,k_2)$-events.

\end{abstract}

\vspace*{.5cm} \noindent {\bf{Key words:}} \emph{\small Poisson distribution, compound Poisson distribution, translated Poisson
distribution, negative binomial distribution, binomial
distribution, m-dependent variables, Wasserstein  norm,
non-uniform estimates.}

\vspace*{.5cm} \noindent {\small {\it MSC 2000 Subject
Classification}:
Primary 60F05.   % Central limit and other weak theorems
Secondary 60G50;     % Sums of independent random variables; random walks
}

\newpage
%%%%%%%%%%%%%%%%%%%%%%%%%%%%%%%%%%%%%%%%%%%%%%%%%%%%%%%%%%%%%%%%%%%%%%%%%%%%%%%%%%%%%%%%%%%%%%%%%%%%%%%%%%%%%%%%%%
%%%%%%%%%%%%%%%%%%%%%%%%%%%%%%%%%%%%%%%%%%%%%%%%%%%%%%%%%%%%%%%%%%%%%%%
%%%% INTRO
%%%%%%%%%%%%%%%%%%%%%%%%%%%%%%%%%%%%%%%%%%%%%%%%%%%%%%%%%%%%%%%%%%%%%%%%%

\section{Introduction} % Initial capital letter, then lower case. No full stop.

Nonuniform estimates for normal approximation are well known, see
the classical results in Chapter 5 of \citep{petrov95} and the references
 \citep{CHSH01}, \citep{CHSH04} and \citep{NefShev12} for some recent developments.
On the other hand, nonuniform estimates for discrete
approximations are only a few. For example, the Poisson approximation to Poisson binomial
distribution has been considered in \cite{Nea03b} and translated Poisson approximation
for independent lattice summands via the Stein method has been discussed in \citep{BCH04}.
 Some general estimates for independent summands under assumption of
matching of pseudomoments were obtained in \citep{ce93}. For
possibly dependent Bernoulli variables, nonuniform estimates for Poisson
approximation problems  were discussed in \citep{TeeSa07}. However, the
estimates obtained had a better accuracy than estimates in total
variation only for $x$ larger than exponent of the sum's mean.
In \citep{cepe11}, 2-runs statistic was approximated by compound
Poisson distribution. In this paper, we obtain nonuniform estimates
for Poisson, compound Poisson, translated Poisson, negative
binomial and binomial approximations, under a quite general set of
assumptions.

We recall that the sequence of random variables $\{X_k \}_{k \geq
1}$ is called $m$-dependent if, for $1 < s < t < \infty$, $t- s > m$, the sigma algebras
generated by $X_1,\dots,X_s$ and $X_t,X_{t+1}\dots$ are
independent. Without loss of generality, we can reduce the sum of $m$-dependent variables to
the sum of 1-dependent ones, by grouping consecutive
$m$ summands. Therefore,  we  consider henceforth, without loss of generality, the sum
$S_n= X_1+X_2+\cdots+X_n$ of non-identically distributed
1-dependent random variables concentrated on nonnegative integers.

We denote the distribution function and the characteristic
function of $S_n$ by $F_n (x)$ and $\w F_n(t)$, respectively.
Similarly, for a signed measure $M$ concentrated on the set $\NN$
of nonnegative integers, we denote by $M(x)=\sum_{k=0}^xM\{k\}$
and $\w M(t)=\sum_{k=0}^\infty \ee^{\ii t k}M\{k\}$, the analogues
of distribution function and Fourier-Stieltjes transform,
respectively. Though our aim is to obtain the non-uniform
estimates, we obtain also estimates for Wasserstein norm defined
as
\[
\norm{M}_W=\sum_{j=0}^\infty\ab{M(j)}.
\]
Note that Wasserstein norm is stronger than total variation norm
defined by $\norm{M}=\sum_{j=0}^\infty\ab{M\{j\}}$.

Next we introduce the approximations considered in this paper.
%For our purposes,
%it is more convenient to use factorial cumulants, rather than
%moments.
Let
\begin{equation*}
\lambda=\Expect S_n,\quad \Gamma_2=\frac{1}{2}(\Var S_n-\Expect
S_n).
\end{equation*}
For brevity, let $z(t)=\ee^{\ii t}-1$. Also, let $\Pi$  and  $\Pi_1$ respectively denote
 the Poisson
distribution with parameter $\lambda$ and its second order
difference multiplied by $\Gamma_2$.  More precisely,
\[\w\Pi(t)=\exponent{\lambda z},\quad\w\Pi_1(t)=\w\Pi(t)\Gamma_2z^2.\]
It is clear that $\Pi+\Pi_1$ is  second-order  (and, consequently,
two-parametric) Poisson approximation. As an alternative to the Poisson
based two-parametric approximation, we choose
 compound Poisson measure $\G$ with the following
Fourier-Stieltjes transform
\begin{equation*}
\w \G(t)=\exponent{\lambda z+\Gamma_2z^2}.
\end{equation*}
 The approximation $\G$ was used in
 many  papers, see \citep{BC02}, \citep{BaXi99}, \citep{Roo03}  and the references therein.
If $\Gamma_2<0$, then $\G$ becomes signed measure, which is not
always convenient and natural for approximation to nonnegative
$S_n$. Therefore,  we define next three distributional
approximations. Translated Poisson ($TP$) approximation  has the
following characteristic function:
\begin{equation*}
\w\TP(t)=\exponent{\lfloor -2\Gamma_2\rfloor \ii t+
(\lambda+2\Gamma_2+\tilde\delta)z}=\exponent{\lambda z+
(2\Gamma_2+\tilde\delta)(z-\ii t)}.
\end{equation*}
Here $\lfloor -2\Gamma_2\rfloor$ and $\tilde\delta$ are respectively the integer
part and the fractional part of $-2\Gamma_2$,  so that
$-2\Gamma_2= \lfloor -2\Gamma_2\rfloor  +\tilde\delta, \quad 0\leqslant
\tilde\delta<1$. The TP approximation was investigated
in numerous papers, see, for example, \citep{BC02}, \citep{BCH04},
\citep{Rol05} and \citep{Rol07}. If
$\Expect S_n<\Var S_n$, then one can apply the negative binomial
approximation, which is defined in the following way:
\begin{equation*}\label{nbinom}
\NB\{j\}=\frac{\Gamma(r+j)}{j!\Gamma(r)}\,\qubar^r(1-\qubar)^j,\quad
(j\in\Zpl),\qquad \frac{r(1-\qubar)}{\qubar}=\lambda,\quad
r\bigg(\frac{1-\qubar}{\qubar}\bigg)^2=2\Gamma_2.
\end{equation*}
Note that
\[
\w\NB(t)=\bigg(\frac{\qubar}{1-(1-\qubar)\ee^{\ii t}}
\bigg)^r=\bigg(1-\frac{(1-\qubar)z}{\qubar}\bigg)^{-r}.
\]
If $\Var S_n<\Expect S_n$, the more natural approximation is the
binomial one defined as follows:
\begin{equation*}\label{binom}
\w\BI(t)=(1+\pbar z)^N,\quad N=\floor{\tilde N}, \quad \tilde N
=\frac{\lambda^2}{2\ab{\Gamma_2}}, \quad\bar{p}=\frac{\lambda}{N}.
\end{equation*}
  Note that  symbols $\qubar$ and $\pbar$ are not related and, in
general, $\qubar+\pbar\neq 1$.

Finally, we introduce some technical notations, related to the
method of proof. Let $\{Y_{k}\}_{k \geq 1}$ be a sequence of
arbitrary real or complex-valued random variables. We assume that
$\w\Expect (Y_1)=\Expect Y_1$ and, for $k\geqslant 2$, define
$\w\Expect (Y_1, Y_2, \cdots, Y_k)$  by
\begin{equation*}
 \w\Expect (Y_1,Y_2,\cdots, Y_k)=\Expect Y_1Y_2\cdots
Y_k-\sum_{j=1}^{k-1}\w\Expect (Y_1,\cdots ,Y_j)\Expect
Y_{j+1}\cdots Y_{k}. \label{capY}
\end{equation*}
 Let
\begin{eqnarray*}
\wE^{+}(X_1)&=&\Expect X_1,\qquad\wE^{+}(X_1,X_2)=\Expect
X_1X_2+\Expect X_1\Expect X_2,\nonumber\\
 \wE^{+}(X_1,\dots,X_k)&=&\Expect X_1\dots
X_k+\sum_{j=1}^{k-1}\wE^{+}(X_1, \dots,X_j)\Expect
X_{j+1}X_{j+2}\cdots X_k,\\
\wE^{+}_2(X_{k-1},X_k)&=&\wexp(X_{k-1}(X_{k-1}-1),X_k)+\wexp(X_{k-1},X_k(X_k-1)),\\
\wexp_2(X_{k-2},X_{k-1},X_k)&=&\wexp(X_{k-2}(X_{k-2}-1),X_{k-1},X_k)+\wexp(X_{k-2},X_{k-1}(X_{k-1}-1),X_k).
\end{eqnarray*}
 We define $j$-th factorial moment of $X_k$ by
$\nu_j(k)=\Expect X_k(X_k-1)\cdots(X_k-j+1)$, ($k=1,2,\dots,n$,
$j=1,2,\dots$). For the sake of convenience, we assume that
$X_k\equiv 0$ and $\nu_j(k)=0$ if $k\leqslant 0$ and $\sum_k^n=0$
if $k>n$. Next we define  remainder terms $R_0$ and $R_1$, which appear in the main results, as
\begin{eqnarray*}
R_0&=&\sum_{k=1}^n\Big\{\nu_2(k)+\nu_1^2(k)+\Expect X_{k-1}X_k\Big\},\\
R_1&=&\sum_{k=1}^n\Big\{\nu_1^3(k)+\nu_1(k)\nu_2(k)+\nu_3(k)+[\nu_1(k-2)+\nu_1(k-1)+\nu_1(k)]\Expect
X_{k-1}X_k\nonumber\\
&&+\wE^{+}_2(X_{k-1},X_k)+\wE^{+}(X_{k-2},X_{k-1},X_k)\Big\}.
\end{eqnarray*}
We use symbol $C$ to denote (in general different) positive
absolute constants.

%%%%%%%%%%%%%%%%%%%%%%%%%%%%%%%%%%%%%%%%%%%%%%%%%%%%%%%%%%%%%%%%%%%%%%%%%%%%%%%%%%%%%%%%%%%%%%%%%%%%%%%%%
%%%%%%%%%%%%%%%%%%%%%%%%%%%%%%%%%%%%%%%%%%%%%%%%%%%%%%%%%%%%%%%%%%%%%%%%%%%%%%%%%%%%%%%%%%%%%%%%%%%%%%%%%
%%%%%%%%%%%%%%%%%%%%%%%%%%%%%%%%%%%%%%%%%%%%%%%%%%%%%%%%%%%%%%%%%%%%%%%%%%%%%%%%%%%%%%%%%%%%%%%%%%%%%%%%%
%%RESULTS
%%%%%%%%%%%%%%%%%%%%%%%%%%%%%%%%%%%%%%%%%%%%%%%%%%%%%%%%%%%%%%%%%%%%%%%%%%%%%%%%%%%%%%%%%%%%%%%%%%%%%%%%%
%%%%%%%%%%%%%%%%%%%%%%%%%%%%%%%%%%%%%%%%%%%%%%%%%%%%%%%%%%%%%%%%%%%%%%%%%%%%%%%%%%%%%%%%%%%%%%%%%%%%%%%%%

\section{The Main Results}

 All the results are obtained under the following conditions:
\begin{eqnarray}
\nu_1(k)&\leqslant& 1/100, \quad
\nu_2(k)\leqslant \nu_1(k),\quad \ab{X_k}\leqslant C_0,\quad (k=1,2,\dots,n),\label{nu12}\\
\lambda&\geqslant&1,\qquad\sum_{k=1}^n\nu_2(k)\leqslant\frac{\lambda}{20},\qquad
\sum_{k=2}^n\ab{Cov(X_{k-1},X_k)}\leqslant\frac{\lambda}{20}.\label{3ab}
 \end{eqnarray}

Assumptions (\ref{nu12}) and (\ref{3ab}) are rather restrictive.
However, they (a) allow to include  independent random variables
as  partial case of general results and (b) are satisfied for many
cases of $k$-runs and $(k_1,k_2)$ events. The method of proof does
not allow to get small constants. Therefore, we have concentrated our
efforts on the order of the accuracy of approximation.
 Next, we state the main results of this paper.
%%%%%%%%%%%%%%%%%%%%%%%%%%%%%%%%%%%%%%%%%%%%%%%%%%%%%%%%%%%%%%%%%%%%%%%%%%%%%%%%%%%%%%%%%%%%%%%%
%%% NONUNIFORM THEOREM
%%%%%%%%%%%%%%%%%%%%%%%%%%%%%%%%%%%%%%%%%%%%%%%%%%%%%%%%%%%%%%%%%%%%%%%%%%%%%%%%%%%%%%%%%%%%%%%%

\begin{theorem} \label{NON} Let conditions (\ref{nu12}) and (\ref{3ab}) be
satisfied. Then, for any $x\in\NN$,
\begin{eqnarray}
\bigg(1+\frac{(x-\lambda)^2}{\lambda}\bigg)\ab{F_n(x)-\Pi(x)}&\leqslant&
C_1\frac{R_0}{\lambda},\label{pi}\\
\bigg(1+\frac{(x-\lambda)^2}{\lambda}\bigg)\ab{F_n(x)-\Pi(x)-\Pi_1(x)}&\leqslant&
C_2\bigg(\frac{R_0^2}{\lambda^{2}}+\frac{R_1}{\lambda\sqrt{\lambda}}\bigg),\label{pi2}\\
\bigg(1+\frac{(x-\lambda)^2}{\lambda}\bigg)\ab{F_n(x)-\G(x)}&\leqslant&
C_3\frac{R_1}{\lambda\sqrt{\lambda}},\label{g}\\
\bigg(1+\frac{(x-\lambda)^2}{\lambda}\bigg)\ab{F_n(x)-\TP(x)}&\leqslant&
C_4\bigg(\frac{R_1+\ab{\Gamma_2}}{\lambda\sqrt{\lambda}}+\frac{\tilde\delta}{\lambda}\bigg).\label{tp}
\end{eqnarray}
If in addition $\Gamma_2>0$, then
\begin{equation}
\bigg(1+\frac{(x-\lambda)^2}{\lambda}\bigg)\ab{F_n(x)-\NB(x)}\leqslant
C_5\bigg(\frac{R_1}{\lambda\sqrt{\lambda}}+\frac{\Gamma_2^2}{\lambda^2\sqrt{\lambda}}\bigg).\label{nb}
\end{equation}
If instead $\Gamma_2<0$, then
\begin{equation}
\bigg(1+\frac{(x-\lambda)^2}{\lambda}\bigg)\ab{F_n(x)-\BI(x)}\leqslant
C_6\bigg(\frac{R_1}{\lambda\sqrt{\lambda}}+\frac{\Gamma_2^2}{\lambda^2\sqrt{\lambda}}\bigg).\label{bi}
\end{equation}
\end{theorem}

\begin{remark}
Nonuniform normal estimates usually match estimates in Kolmogorov
metric. Similarly, the bounds in (\ref{pi})-(\ref{bi})  match estimates in
total variation:
\[\norm{F_n-\Pi}\leqslant
C_7\frac{R_0}{\lambda},\quad \norm{F_n-\Pi-\Pi_1}\leqslant
C_8\bigg(\frac{R_0^2}{\lambda^{2}}+\frac{R_1}{\lambda\sqrt{\lambda}}\bigg),\quad
\norm{F_n-\G}\leqslant C_9\frac{R_1}{\lambda\sqrt{\lambda}},
\]
and etc., see \citep{CeVe13}.
\end{remark}

Estimates for Wasserstein metric easily follow by summing up
nonuniform estimates.

%%%%%%%%%%%%%%%%%%%%%%%%%%%%%%%%%%%%%%%%%%%%%%%%%%%%%%%%%%%%%%%%%%%%%%%%%%%%%%%%%%%%%%%%%%%%%%%%
%%% WASSER THEOREM
%%%%%%%%%%%%%%%%%%%%%%%%%%%%%%%%%%%%%%%%%%%%%%%%%%%%%%%%%%%%%%%%%%%%%%%%%%%%%%%%%%%%%%%%%%%%%%%%

\begin{theorem} \label{WASS} Let conditions (\ref{nu12}) and (\ref{3ab}) be
satisfied. Then,
\begin{eqnarray}
\norm{F_n-\Pi}_W&\leqslant&
C_{10}\frac{R_0}{\sqrt{\lambda}},\label{wpi}\\
\norm{F_n-\Pi-\Pi_1}_W&\leqslant&
C_{11}\bigg(\frac{R_0^2}{\lambda\sqrt{\lambda}}+\frac{R_1}{\lambda}\bigg),\label{wpi2}\\
\norm{F_n-\G}_W&\leqslant&
C_{12}\frac{R_1}{\lambda},\label{wg}\\
\norm{F_n-\TP}_W&\leqslant&
C_{13}\bigg(\frac{R_1+\ab{\Gamma_2}}{\lambda}+\frac{\tilde\delta}{\sqrt{\lambda}}\bigg).\label{wtp}
\end{eqnarray}
When in addition $\Gamma_2>0$, we have
\begin{equation}
\norm{F_n-\NB}_W\leqslant
C_{14}\bigg(\frac{R_1}{\lambda}+\frac{\Gamma_2^2}{\lambda^2}\bigg),\label{wnb}
\end{equation}
and when $\Gamma_2<0$, we have
\begin{equation}
\norm{F_n-\BI}_W\leqslant
C_{15}\bigg(\frac{R_1}{\lambda}+\frac{\Gamma_2^2}{\lambda^2}\bigg).\label{wbi}
\end{equation}
\end{theorem}

Observe that the local nonuniform estimates have better order of accuracy.
%%%%%%%%%%%%%%%%%%%%%%%%%%%%%%%%%%%%%%%%%%%%%%%%%%%%%%%%%%%%%%%%%%%%%%%%%%%%%%%%%%%%%%%%%%%%%%%%
%%% LOCAL NONUNIFORM THEOREM
%%%%%%%%%%%%%%%%%%%%%%%%%%%%%%%%%%%%%%%%%%%%%%%%%%%%%%%%%%%%%%%%%%%%%%%%%%%%%%%%%%%%%%%%%%%%%%%%

\begin{theorem} \label{NONL} Let conditions (\ref{nu12}) and (\ref{3ab}) hold. Then, for any $x\in\NN$,
\begin{eqnarray}
\bigg(1+\frac{(x-\lambda)^2}{\lambda}\bigg)\ab{F_n\{x\}-\Pi\{x\}}&\leqslant&
C_{16}\frac{R_0}{\lambda\sqrt{\lambda}},\label{pil}\\
\bigg(1+\frac{(x-\lambda)^2}{\lambda}\bigg)\ab{F_n\{x\}-\Pi\{x\}-\Pi_1\{x\}}&\leqslant&
C_{17}\bigg(\frac{R_0^2}{\lambda^{2}\sqrt{\lambda}}+\frac{R_1}{\lambda^2}\bigg),\label{pi2l}\\
\bigg(1+\frac{(x-\lambda)^2}{\lambda}\bigg)\ab{F_n\{x\}-\G\{x\}}&\leqslant&
C_{18}\frac{R_1}{\lambda^2},\label{gl}\\
\bigg(1+\frac{(x-\lambda)^2}{\lambda}\bigg)\ab{F_n\{x\}-\TP\{x\}}&\leqslant&
C_{19}\bigg(\frac{R_1+\ab{\Gamma_2}}{\lambda^2}+\frac{\tilde\delta}{\lambda\sqrt{\lambda}}\bigg).\label{tpl}
\end{eqnarray}
If in addition $\Gamma_2>0$, then
\begin{equation}
\bigg(1+\frac{(x-\lambda)^2}{\lambda}\bigg)\ab{F_n\{x\}-\NB\{x\}}\leqslant
C_{20}\bigg(\frac{R_1}{\lambda^2}+\frac{\Gamma_2^2}{\lambda^3}\bigg).\label{nbl}
\end{equation}
If instead $\Gamma_2<0$, then
\begin{equation}
\bigg(1+\frac{(x-\lambda)^2}{\lambda}\bigg)\ab{F_n\{x\}-\BI\{x\}}\leqslant
C_{21}\bigg(\frac{R_1}{\lambda^2}+\frac{\Gamma_2^2}{\lambda^3}\bigg).\label{bil}
\end{equation}
\end{theorem}

\begin{remark} (i) Estimates in (\ref{pil})-(\ref{bil}) match estimates in local
metric, see \citep{CeVe13}.

\noindent (ii) Consider the case of independent Bernoulli
variables with $p\leqslant 1/20$ and $\lambda\geqslant 1$. Then,
for all integers $x$, Poisson approximation is of the order
\[
\frac{C\sum_{j=1}^np_j^2}{(1+(x-\lambda)^2/\lambda)\lambda\sqrt{\lambda}},\]
which is usually much better than
\[
\min(x^{-1},\lambda^{-1})\sum_{j=1}^np_j^2
\]
from \citep{Nea03a}.
\end{remark}
%%%%%%%%%%%%%%%%%%%%%%%%%%%%%%%%%%%%%%%%%%%%%%%%%%%%%%%%%%%%%%%%%%%%%%%%%%%%%%%%%%%%%%%%%%%%%%%%

%%%%%%%%%%%%%%%%%%%%%%%%%%%%%%%%%%%%%%%%%%%%%%%%%%%%%%%%%%%%%%%%%%%%%%%%%%%%%%%%%%%%%%%%%%%%%%%%%%%%%
\section{Some Applications}

\textbf{ (i): Asymptotically sharp constant for Poisson approximation
to Poisson binomial distribution.} Formally, independent random
variables make a subset of 1-dependent variables. Therefore, one
can rightly expect that results of the previous section apply to
independent summands as well. We exemplify this fact by considering one of
the best known cases in Poisson approximation theory. Let
$W=\xi_1+\xi_2+\cdots+\xi_n$, where $\xi_i$ are independent
Bernoulli variables with $P(\xi_i=1)=1-P(\xi_i=0)=p_i$. Let
$\lambda=\sum_{1}^np_i$, $\lambda_2=\sum_{1}^np_i^2$. As shown in
\citep{BaXi06} (see equation (1.8)),
\begin{equation}
\norm{\eL(W)-\Pi}_W\leqslant\frac{1.1437\lambda_2}{\sqrt{\lambda}}.
\label{B}
\end{equation}
Though absolute constant in (\ref{B}) is small, we shall show that
asymptotically sharp constant is much smaller. Let $\max\limits_{i}p_i\to
0$ and $\lambda\to\infty$, as $n\to\infty$. Then
\begin{equation}
\lim_{n\to\infty}\frac{\sqrt{\lambda}}{\lambda_2}\norm{\eL(W)-\Pi}_W=\frac{1}{\sqrt{2\pi}}\leqslant
0.399.\label{sharp}
\end{equation}
Indeed, we have
\[\Ab{\norm{\eL(W)-\Pi}_W-\frac{\lambda_2}{\sqrt{2\pi\lambda }}}
\leqslant\norm{\eL(W)-\Pi-\Pi_1}_W+\Ab{\norm{\Pi_1}_W-\frac{\lambda_2}{\sqrt{2\pi\lambda
}}}.\]

If $\max\limits_{i} p_i\leqslant 1/20$ and $\lambda\geqslant 1$, then it follows from
(\ref{wpi2}) that
\[
\norm{\eL(W)-\Pi-\Pi_1}_W\leqslant
\frac{C\lambda_2}{\sqrt{\lambda}}\bigg(\max_jp_j+\frac{1}{\sqrt{\lambda}}\bigg).
\]
For the estimation of the second difference, we require some
notations for measures. Let $Z$ be a measure, corresponding to Fourier-Stieltjes
transform $z(t)= (e^{it}-1)$. Let  product and powers of measures be
understood in the convolution sense. Then, by the properties of
norms and Proposition 4 from \citep{Roo99}  (see also Lemma 6.2 in
\citep{CeVe13})), we get
\begin{eqnarray*}
\Ab{\norm{\Pi_1}_W-\frac{\lambda_2}{\sqrt{2\pi\lambda}}}&=&
\Ab{\frac{\lambda_2}{2}\norm{\Pi
Z^2}_W-\frac{\lambda_2}{\sqrt{2\pi\lambda
}}}=\frac{\lambda_2}{2}\Ab{\norm{\Pi
Z^2}_W-\frac{\sqrt{2/\pi}}{\sqrt{\lambda}}}\\
&=& \frac{\lambda_2}{2}\Ab{\norm{\Pi
Z}-\frac{\sqrt{2/\pi}}{\sqrt{\lambda}}}\leqslant
\frac{C\lambda_2}{2\lambda}=\frac{\lambda_2}{\sqrt{\lambda}}\frac{C}{2\sqrt{\lambda}}.
\end{eqnarray*}
Thus, for $\max\limits_{i}p_i\leqslant 1/20$ and  $\lambda\geqslant 1$, we
obtain asymptotically sharp norm estimate
\begin{equation*}
\Ab{\norm{\eL(W)-\Pi}_W-\frac{\lambda_2}{\sqrt{2\pi\lambda
}}}\leqslant
\frac{C\lambda_2}{\sqrt{\lambda}}\bigg(\max_jp_j+\frac{1}{\sqrt{\lambda}}\bigg),
\label{sharpa}
\end{equation*}
which is even more general than (\ref{sharp}).

%%%%%%%%%%%%%%%%%%%%%%%%%%%%%%%%%%%%%%%%%%%%%%%%%%%%%%%%%%%%%%%%%%%%%%%%%%%%%%%%%%%%%%%%%%%%%%%%%%%%%%%%%%%%%%%%%%%%%%
%%%%%%%%%%%%%%%%%%%%%%%%%%%%%%%%%%%%%%%%%%%%%%%%%%%%%%%%%%%%%%%%%%%%%%%%%%%%%%%%%%%%%%%%%%%%%%%%%%%%%%%%%%%%%%%%%%%%%%
%%%%%%%%%%%%%%%%%%%%%%%%%%%%%%%%%%%%%%%%%%%%%%%%%%%%%%%%%%%%%%%%%%%%%%%%%%%%%%%%%%%%%%%%%%%%%%%%%%%%%%%%%%%%%%%%%%%%%%

\noindent \textbf{ (ii): Negative binomial approximation to 2-runs.}
The k-runs (and especially 2-runs) statistic is one of the best
investigated cases of sums of dependent discrete random variables,
see \citep{WXia08} and the references therein. Let
$S_n=X_1+X_2+\dots+X_n$, where $X_i=\eta_i\eta_{i+1}$ and
$\eta_j\sim Be(p)$, $(j=1,2,\dots,n+1)$ are independent Bernoulli
variables. Then $S_n$ is called 2-runs statistic. It is known that
then
\[
\lambda=np^2,\quad\Gamma_2= \frac{np^3(2-3p)-2p^3(1-p)}{2}.
\]
Let $p\leqslant 1/20$ and $np^2\geqslant 1$. Then, from (\ref{nb}) it
follows for any $x\in\NN$,
\begin{equation*}
\bigg(1+\frac{(x-\lambda)^2}{\lambda}\bigg)\ab{F_n(x)-\NB(x)}\leqslant
C\frac{p}{\sqrt{n}}. \label{nba}\end{equation*}
  This estimate has
the same order as the estimate in total variation, see
 and \citep{BrX01}  and \citep{CeVe13}.

\noindent \textbf{(iii): Binomial approximation to $(k_1,k_2)$-events.} Let $\eta_i$
be independent Bernoulli $Be(p)$ ($0<p<1$) variables and let
$Y_j=(1-\eta_{j-m+1})\cdots
(1-\eta_{j-k_2})\eta_{j-k_2+1}\cdots\eta_{j-1}\eta_j$,
$j=m,m+1,\dots,n$, $k_1+k_2=m$. Further, we assume that $k_1>0$
and $k_2>0$. Then $N(n;k_1,k_2)=Y_m+Y_{m+1}+\cdots+Y_n$ denote the number of
$(k_1,k_2)$-events and we denote
its distribution  by $\Ha$. The Poisson approximation to $H$ has been considered in \cite{vs04}.  Let $a(p)=(1-p)^{k_1}p^{k_2}$.

Note that $Y_1,Y_2,\dots$ are $m$-dependent.  However, one can
group summands in the following natural way:
\begin{eqnarray*}
N(n;k_1,k_2)&=&(Y_m+Y_{m+1}+\cdots+Y_{2m-1})+(Y_{2m}+Y_{2m+1}+\cdots+Y_{3m-1})+\dots \\
 &=& X_1+X_2+\dots .
\end{eqnarray*}
Each $X_j$, with probable exception of the last one, contains $m$
 summands. It is not difficult  to check that $X_1,X_2,\dots$ are
 1-dependent
 Bernoulli variables.
Then all parameters can be written explicitly. Set $N=\lfloor
\tilde N\rfloor$ be the integer part of $\tilde N$ defined by
\[ \tilde
N=\frac{(n-m+1)^2}{(n-m+1)(2m-1)-m(m-1)},\quad\bar{p}=\frac{(n-m+1)a(p)}{N}.\]
 It is
known (see \citep{CeVe13}) that
 \[
  \lambda=(n-m+1)a(p),\quad
 \Gamma_2=-\frac{a^2(p)}{2}[(n-m+1)(2m-1)-m(m-1)],
 \quad R_1\leqslant
C(n-m+1)m^2a^3(p).
\]

\noindent  Let now $ \lambda \geq 1$ and
$ma(p)\leqslant 0.01$. Then  it follows from (\ref{bi}) that, for
any $x\in\NN$,
\begin{equation*}
\bigg(1+\frac{(x-\lambda)^2}{\lambda}\bigg)\ab{\Ha(x)-\BI(x)}\leqslant
C \frac{a^{3/2}(p)m^2}{\sqrt{n-m+1}}. \label{bia}\end{equation*}

%%%%%%%%%%%%%%%%%%%%%%%%%%%%%%%%%%%%%%%%%%%%%%%%%%%%%%%%%%%%%%%%%%%%%%%%%%%%%%%%%%%%%%%%%%%%%%%%%%%%%%%%%%%%%%%%%%%%%
%%%%%%%%%%%%%%%%%%%%%%%%%%%%%%%%%%%%%%%%%%%%%%%%%%%%%%%%%%%%%%%%%%%%%%%%%%%%%%%%%%%%%%%%%%%%%%%%%%%%%%%%%%%%%%%%%%%%%
% AUXILARY
%%%%%%%%%%%%%%%%%%%%%%%%%%%%%%%%%%%%%%%%%%%%%%%%%%%%%%%%%%%%%%%%%%%%%%%%%%%%%%%%%%%%%%%%%%%%%%%%%%%%%%%%%%%%%%%%%%%%%
%%%%%%%%%%%%%%%%%%%%%%%%%%%%%%%%%%%%%%%%%%%%%%%%%%%%%%%%%%%%%%%%%%%%%%%%%%%%%%%%%%%%%%%%%%%%%%%%%%%%%%%%%%%%%%%%%%%%%
%%%%%%%%%%%%%%%%%%%%%%%%%%%%%%%%%%%%%%%%%%%%%%%%%%%%%%%%
\section{Auxiliary results}
Let $\theta$  to denote a real or complex quantity
satisfying $\ab{\theta}\leqslant 1$.  Moreover, let
$Z_j=\exponent{\ii t X_j}-1$, $\Psi_{j,k}=\wE(Z_j,Z_{j+1},\dots,Z_k)$.
As before, we assume that $\nu_j(k)=0$ and $X_k=0$ for $k\leqslant
0$ and $z(t)=\eit-1$.
%Moreover, to make our notation
%shorter we write $\wE_z (1,2)$ instead of $\wE(Z_1,Z_2)$,
% $\wE_z(1,2,3)$ instead of $\wE(Z_1,Z_2,Z_3)$ and so on.
 Also, we omit the argument $t$, wherever possible and,  for example,
write $z$ instead of $z(t)$. Hereafter, the primes denote the derivatives with respect to $t$.

%%%%%%%%%%%%%%%%%%%%%%%%%%%%%%%%%%%%%%%%%%%%%%%%%%%%%%%%%%%%%%%%%%%%%%%%%%%%%%%%%%%%%%%%%%%%%%%%%%%%%%%%%%%%%%%%%%%%%%%%
%%%% FACTORIAL
%%%%%%%%%%%%%%%%%%%%%%%%%%%%%%%%%%%%%%%%%%%%%%%%%%%%%%%%%%%%%%%%%%%%%%%%%%%%%
\begin{lemma} \label{factorial} Let $X$ be concentrated on nonnegative integers and $\nu_3<\infty$. Then, for all $t\in\RR$,
\begin{eqnarray*}
\Expect \exponent{\ii t X} &=&
1+\nu_1z+\nu_2\frac{z^2}{2}+\theta\frac{\nu_3\ab{z}^3}{6},\\
\Expect (\exponent{\ii t X})' &=&
\nu_1z'+\nu_2\frac{(z^2)'}{2}+\theta\frac{\nu_3\ab{z}^2}{2},\\
\Expect (\exponent{\ii t X})'' &=&
\nu_1z''+\nu_2\frac{(z^2)''}{2}+\theta 2\nu_3\ab{z}.
\end{eqnarray*}
\end{lemma}

\noindent \Proof First equality is well known expansion of characteristic
function in factorial moments. The other two equalities also
easily follow from expansions in powers of $z$. For example,
\begin{eqnarray}
(\ee^{\ii t X})''&=&\ii^2X^2\ee^{\ii t
X}=\ii^2X(X-1)(\eit)^2\ee^{\ii t(X-2)}+\ii^2\eit X\ee^{\ii
t(X-1)}\nonumber\\
&=&\ii^2(\eit)^2X(X-1)[1+\theta(X-2)\ab{z}]+\ii^2\eit
X[1+(X-1)z+\theta
(X-1)(X-2)\ab{z}^2/2]\nonumber\\
&=&Xz''+\frac{X(X-1)}{2}(z^2)''+\theta
2X(X-1)(X-2)\ab{z}.   \rlap{$\qquad \Box$} \label{dui}
\end{eqnarray}

%%%%%%%%%%%%%%%%%%%%%%%%%%%%%%%%%%%%%%%%%%%%%%%%%%%%%%%%%%%%%%%%%%%%%%%%%
%%%%\WE HOLDER
%%%%%%%%%%%%%%%%%%%%%%%%%%%%%%%%%%%%%%%%%%%%%%%%%%%%%%%%%%%%%%%%%%%%%%%%%%%%%%%%%%%%%%%%%%%%%%%%%%%%%%%%%%%%%%%%%%
\begin{lemma}(\citep{H82}) \label{Hei3aa} Let $Y_1,Y_2,\dots,Y_k$ be 1-dependent complex-valued random variables with
$\Expect\ab{Y_m}^2<\infty,~ 1 \leq m \leq k. $ Then
\[
\ab{\wE (Y_1, Y_2, \cdots, Y_k)}\leqslant
2^{k-1}\prod_{m=1}^k(\Expect\ab{Y_m}^2)^{1/2}.
\]
\end{lemma}
%%%%%%%%%%%%%%%%%%%%%%%%%%%%%%%%%%%%%%%%%%%%%%%%%%%%%%%%%%%%%%%%%%%%%%
%%%%%%%%%%%%%%%%%%%%%%%%

%\WE ISV ZJ
%%%%%%%%%%%%%%%%%%%%%%%%%%%%%%%%%%%%%%%%%%%%%%%%%%%%%%%%%%%%%%%%%%%%%%%%%%%%%%%%%%%%%%%%%%%%%%%%%%%%%%%%%%%%%%%%%%%%%%
\begin{lemma} Let conditions (\ref{nu12})  be
satisfied and $j<k-1$. Then, for all $t$,
\begin{eqnarray}
\ab{\Psi_{j,k}}&\leqslant&
4^{k-j}\ab{z}\prod_{l=j}^k\sqrt{\nu_1(l)},\label{psi0}\\
\ab{\Psi_{j,k}'}&\leqslant&
4^{k-j}\ab{z}(k-j+1)\prod_{l=j}^k\sqrt{\nu_1(l)},\label{psi1}\\
\ab{\Psi_{j,k}''}&\leqslant&
\sqrt{2}C_04^{k-j}\ab{z}(k-j+1)(k-j)\prod_{l=j}^k\sqrt{\nu_1(l)}.\label{psi2}
\end{eqnarray}
\end{lemma}
\noindent \Proof First two estimates follow from more general
estimates  in (47) and Lemma 7.5 in \citep{CeVe13}. Note also the following inequalities:
\begin{equation}
\ab{z}\leqslant 2,\quad\ab{Z_k}\leqslant 2,\quad \ab{Z_k}\leqslant
X_k\ab{z},\quad \Expect X_i^2=\nu_2(i)+\nu_1(i)\leqslant
2\nu_1(i). \label{ntrivial}
\end{equation}
Therefore, by Lemma \ref{Hei3aa} and for $m \leq k$,
\begin{eqnarray*}
\lefteqn{\ab{\wE(Z_j,\dots,Z_m',\dots,Z_i',\dots Z_k)}\leqslant
2^{k-j}\sqrt{\Expect\ab{Z_m'}^2\Expect\ab{Z_i'}^2}\prod_{l=j,l\ne
m,i}^k\sqrt{2\ab{z}\nu_1(l)}}\hskip 0.5cm\\
&\leqslant&2^{k-j}\sqrt{2\nu_1(m)2\nu_1(i)}2^{(k-j-1)/2}\ab{z}^{(k-j-1)/2}
\prod_{l=j,l\ne m,i}^k\sqrt{\nu_1(l)}\leqslant
4^{k-j}2^{-1}\ab{z}\prod_{l=j}^k\sqrt{\nu_1(l)}.
\end{eqnarray*}
Similarly,
\begin{eqnarray*}
\ab{\wE(Z_j,\dots,Z_i'',\dots, Z_k)}&\leqslant&
2^{k-j}\sqrt{\Expect\ab{Z_i''}^2}\prod_{l=j,l\ne
i}^k\sqrt{2\ab{z}\nu_1(l)}\\
&\leqslant& 2^{k-j}\sqrt{\Expect
X_i^4}2^{(k-j)/2}\ab{z}^{(k-j)/2}\prod_{l=j,l\ne
i}^k\sqrt{\nu_1(l)}\\
&\leqslant&4^{k-j}2^{-1}\ab{z}C_0\sqrt{\Expect
X_i^2}\prod_{l=j,l\ne i}^k\sqrt{\nu_1(l)}\leqslant
4^{k-j}2^{-1/2}C_0\prod_{l=j}^k\sqrt{\nu_1(l)}.
\end{eqnarray*}
Thus,
\begin{eqnarray*}
\ab{\Psi_{j,k}''}&\leqslant&\sum_{i=j}^k\ab{\wE(Z_j,\dots,Z_i'',\dots,Z_k)}+
\sum_{i=j}^k\sum_{m=j,m\ne
i}^k\ab{\wE(Z_j,\dots,Z_m',\dots,Z_i',\dots,Z_k)}\\
&\leqslant&(k-j+1)4^{k-j}C_02^{-1/2}\ab{z}\prod_{l=j}^k\sqrt{\nu_1(l)}+
(k-j+1)(k-j)4^{k-j}2^{-1}\ab{z}\prod_{l=j}^k\sqrt{\nu_1(l)}\\
&\leqslant&
\sqrt{2}C_04^{k-j}(k-j+1)(k-j)\ab{z}\prod_{l=j}^k\sqrt{\nu_1(l)}. \rlap{$\qquad \Box$}
\end{eqnarray*}

%%%%%%%%%%%%%%%%%%%%%%%%%%%%%%%%%%%%%%%%%%%%%%%%%%%%%%%%%%%%%%%%%%%%%%%%%%%%%%%%%%%%%%%%%%%%%%%%%%%%%%%%%%%%
%%%%%phi old
%%%%%%%%%%%%%%%%%%%%%%%%%%%%%%%%%%%%%%%%%%%%%%%%%%%%%%%%%%%%%%%%%%%%%%%%%%%%%%%%%%%%%%%%%%%%%%%%%%%
%%%%%%%%%%%%%%%%%%%%%%%%%%%%%%%%%%%%%%%%%%%%%%%%%%%%%%%%%%%%%%%%%%%%%%%%%%%%%%%%%%%%%%%%%%%%%%%%%%%
In the following Lemmas \ref{hei}--\ref{fi13}, we  present some facts about
characteristic function $\w F_n(t)$ from \citep{CeVe13}. Here again we
assume (\ref{nu12}), though many relations hold also under weaker
assumptions, see \citep{CeVe13}. We begin from Heinrich's
representation of $\w F_n$ as product of functions.
 %%%%%%%%%%%%%%%%%%%%%%%%%%%%%%%%%%%%%%%%%%%%%%%%%%%%%%%%%%%%%%%%%%%%%%%%%%%%%%%%%%%%%%%%%%%%%%%%%%%%%
%%% HEINRICH
%%%%%%%%%%%%%%%%%%%%%%%%%%%%%%%%%%%%%%%%%%%%%%%%%%%%%%%%%%%%%
\begin{lemma}\label{hei} Let (\ref{nu12}) hold.  Then
 $\w F_n(t)=\vfi_1(t)\vfi_2(t)\dots\vfi_n(t)$,
 where $\vfi_1(t)=\Expect e^{\ii t X_1}$ and, for $k=2,\dots, n$,
\begin{equation}\vfi_k=1+ \Expect
Z_k+\sum_{j=1}^{k-1}\frac{\Psi_{j,k}}{\vfi_j\vfi_{j+1}\dots \vfi_{k-1}}.
\label{hein}\end{equation}
\end{lemma}

%%%%%%%%%%%%%%%%%%%%%%%%%%%%%%%%%%%%%%%%%%%%%%%%%%%%%%%%%%%%%%%%%%%%%%%%%%%%%%%%%%%%%
%%%%LEMA FROM 2013
%%%%%%%%%%%%%%%%%%%%%%%%%%%%%%%%%%%%%%%%%%%%%%%%%

 Let
\begin{eqnarray*}
 g_j(t)&=&\Exponent{\nu_1(j)z+\Big(\frac{\nu_2(j)-\nu_1^2(j)}{2}+\wE(X_{j-1},X_j)
\Big)z^2},\label{gj}\\
%\tilde\vfi_k&=&\vfi_k\exponent{-\ii\nu_1(k)t}, \quad\tilde g_k=
%g_k\exponent{-\ii\nu_1(k)t},\\
\lambda_k&=&1.6\nu_1(k)-0.3\nu_1(k-1)-2\nu_2(k)-0.1\Expect
 X_{k-2}X_{k-1}-15.58\Expect X_{k-1}X_{k},\\
\gamma_2(k)&=&\frac{\nu_2(k)}{2}+\wE(X_{k-1},X_k),\nonumber\\
%r_0(k)&=&\nu_2(k)+\sum_{l=0}^3\nu_1^2(k-l)+\Expect
%X_{k-1}X_k,\label{r0k}\\
r_1(k)&=&\nu_3(k)+\sum_{l=0}^5\nu_1^3(k-l)+[\nu_1(k-1)+\nu_1(k-2)]\Expect
X_{k-1}X_k+\wE_2^{+}(X_{k-1},X_k)\nonumber\\
&&+\wE^{+}(X_{k-2},X_{k-1},X_k),\label{r1kn}
\end{eqnarray*}
\begin{lemma}\label{fi13} Let the conditions in (\ref{nu12})
hold. Then
\begin{eqnarray}
\frac{1}{\ab{\vfi_k}}&\leqslant&\frac{10}{9},\label{sh1}\\
\ab{\vfi_k}&\leqslant&
\exponent{-\lambda_k\sin^2(t/2)},\quad \ab{g_k}\leqslant \exponent{-\lambda_k\sin^2(t/2)}\label{elk}\\
 \frac{1}{\vfi_{k-1}}&=&1+C\theta\ab{z}\{\nu_1(k-2)+\nu_1(k-1)\},\label{f4}\\
\vfi_k'&=&33\theta[\nu_1(k)+\nu_1(k-1)],\label{fis0}\\
\sum_{k=1}^n\ab{\vfi_k-g_k}&\leqslant& CR_1\ab{z}^3,\quad
\sum_{k=1}^n\ab{\vfi_k'-g_k'}\leqslant CR_1\ab{z}^2.\label{sumi}
\end{eqnarray}
\end{lemma}

Similar estimates hold for the second derivative, as seen in the next lemma.

%%%%%%%%%%%%%%%%%%%%%%%%%%%%%%%%%%%%%%%%%%%%%%%%%%%%%%%%%%%%%%%%%%%%%%%%%%%%%%%%%%%%%%%%%%%%%%%%%%%%%%
%%% VFIS2
%%%%%%%%%%%%%%%%%%%%%%%%%%%%%%%%%%%%%%%%%%%%%%%%%%%%%%%%%%%%%%%%%%%%%%%%%%%%%%%%%%%%%%%%%%%%%%%%%%%%%%
\begin{lemma}\label{fin} Let (\ref{nu12})
hold. Then, for $k=1,2,\dots,n$,
\begin{eqnarray}
\vfi_k''&=& \theta C_{22}[\nu_1(k)+\nu_1(k-1)],\label{vfi2}\\
\vfi_k''&=& \nu_1(k)z''+\gamma_2(k)(z^2)''+\theta C\ab{z}r_1(k).
\label{vfi2a}
\end{eqnarray}
\end{lemma}

\Proof From Lemma \ref{hei}, it follows that
\begin{eqnarray}
\vfi_k''&=&(\Expect
Z_k)''+\sum_{j=1}^{k-1}\frac{\Psi_{j,k}''}{\vfi_j\cdots\vfi_{k-1}}-
2\sum_{j=1}^{k-1}\frac{\Psi_{j,k}'}{\vfi_j\cdots\vfi_{k-1}}\sum_{i=j}^{k-1}\frac{\vfi_i'}{\vfi_i}\nonumber\\
&&+\sum_{j=1}^{k-1}\frac{\Psi_{j,k}}{\vfi_j\cdots\vfi_{k-1}}\bigg(\sum_{i=j}^{k-1}\frac{\vfi_i'}{\vfi_i}\bigg)^2+
\sum_{j=1}^{k-1}\frac{\Psi_{j,k}}{\vfi_j\cdots\vfi_{k-1}}\sum_{i=j}^{k-1}\bigg(\frac{\vfi_i'}{\vfi_i}\bigg)^2\nonumber\\
&&-
\sum_{j=1}^{k-1}\frac{\Psi_{j,k}}{\vfi_j\cdots\vfi_{k-1}}\sum_{i=j}^{k-1}\frac{\vfi_i''}{\vfi_i}.
\label{antraisv}
\end{eqnarray}
We prove (\ref{vfi2}) by mathematical induction. Note
that by Lemma \ref{factorial} $(\Expect Z_k)''= C \theta \nu_1(k)$.
Moreover, for $j\leqslant k-2$,
\begin{equation}
\prod_{l=j}^k\sqrt{\nu_1(l)}=\sqrt{\nu_1(k)\nu_1(k-1)}\prod_{l=j}^{k-2}\sqrt{\nu_1(l)}
\leqslant\frac{\nu_1(k)+\nu_1(k-1)}{2}10^{-(k-j-1)}.
\label{pro}
\end{equation}
Applying (\ref{pro}) to (\ref{psi0}), for all $j\leqslant k-2$, we
prove
\begin{equation}
\ab{\Psi_{j,k}}\leqslant
10\bigg(\frac{4}{10}\bigg)^{k-j}[\nu_1(k)+\nu_1(k-1)]. \label{ce1}
\end{equation}
Taking into account (\ref{ntrivial}) and (\ref{nu12}), it is easy
to check that
\begin{eqnarray*}
\ab{\wE(Z_{k-1},Z_k)}&\leqslant& \Expect\ab{Z_{k-1}Z_{k}}+\Expect
\ab{Z_{k-1}}\Expect\ab{Z_k}=
\Expect\ab{Z_{k-1}Z_{k}}/2+\Expect\ab{Z_{k-1}Z_{k}}/2+\Expect \ab{Z_{k-1}}\Expect\ab{Z_k}/2\\
&&+\Expect
\ab{Z_{k-1}}\Expect\ab{Z_k}/2\leqslant\Expect\ab{Z_{k-1}}+\Expect\ab{Z_{k}}+
0.01\Expect\ab{Z_{k-1}}+0.01\Expect\ab{Z_k}\\
&\leqslant& 2.02[\nu_1(k-1)+\nu_1(k)].
\end{eqnarray*}
Therefore, we see that (\ref{ce1}) holds also for $j=k-1$. From
inductional assumption, (\ref{sh1}), (\ref{fis0}) and (\ref{nu12}),
it follows
\[
\frac{\ab{\vfi_i''}}{\ab{\vfi_i}}\leqslant
C_{22}[\nu_1(i-1)+\nu_1(i)]\frac{10}{9}\leqslant
\frac{2C_{22}}{90}.
\]
Using (\ref{sh1}) and the previous estimate, we
obtain
\begin{eqnarray*}
\lefteqn{\Ab{\sum_{j=1}^{k-1}\frac{\Psi_{j,k}}{\vfi_j\cdots\vfi_{k-1}}\sum_{i=j}^{k-1}\frac{\vfi_i''}{\vfi_i}}
\leqslant
\sum_{j=1}^{k-1}\bigg(\frac{10}{9}\bigg)^{k-j}\ab{\Psi_{j,k}}\sum_{i=j}^{k-1}\frac{\ab{\vfi_i''}}{\ab{\vfi_i}}}\hskip
0.5cm\\
&\leqslant& \sum_{j=1}^{k-1}
10\bigg(\frac{4}{9}\bigg)^{k-j}[\nu_1(k)+\nu_1(k-1)](k-j)
\frac{2C_{22}}{90}
\leqslant\frac{8C_{22}}{25}[\nu_1(k)+\nu_1(k-1)].
\end{eqnarray*}
Estimating all other sums (without using induction arguments) in a similar manner, we finally arrive at the estimate
\[\ab{\vfi_k''}\leqslant C_{23}[\nu_1(k-1)+\nu_1(k)]+
\frac{8C_{22}}{25}[\nu_1(k)+\nu_1(k-1)].\]
 It remains to choose $C_{22}=25C_{23}/17$ to complete the proof of
 (\ref{vfi2}).

 \noindent Since the proof of (\ref{vfi2a}) is quite similar, we
 give only a general outline of it. First, we assume that
 $k\geqslant 6$. Then in (\ref{antraisv})  split all sums into
 $\sum_{j=1}^{k-5}+\sum_{j=k-4}^{k-1}$.
 Next, note that
\[
\prod_{l=j}^k\sqrt{\nu_1(l)}\leqslant\prod_{l=k-5}^k\sqrt{\nu_1(l)}\prod_{l=j}^{k-6}\Big(\frac{1}{10}\Big)\leqslant
\sum_{l=k-5}^k\nu_1^3(l)10^{-(k-j-5)}\leqslant
r_1(k)10^{-(k-j-5)}.
\]
Therefore, applying (\ref{psi0})--(\ref{psi2}) and using (\ref{sh1}), (\ref{fis0}) and (\ref{vfi2}),
 we easily prove that all sums
$\sum_{j=1}^{k-5}$ are by absolute value less than
$C\ab{z}r_1(k)$. The cases $j=k-4,k-3,k-2$ all contain at least
three $Z_i$ and can be estimated directly by $C\ab{z}r_1(k)$. For
example,
\[\ab{\wE(Z_{k-3},Z_{k-2},Z_{k-1},Z_k)}\leqslant
4\wE^{+}(\ab{Z_{k-2}},\ab{Z_{k-1}},\ab{Z_k})\leqslant
C\ab{z}^3\wE^{+}(X_{k-2},X_{k-1},X_k)\leqslant C\ab{z}r_1(k).
\]
 Easily verifiable estimates
$\ab{(\wE(Z_{k-1},Z_k))'}\ab{\vfi_{k-1}'}\leqslant C\ab{z}r_1(k)$,
$\ab{\wE(Z_{k-1},Z_k)}\ab{\vfi_{k-1}'}^2\leqslant C\ab{z}r_1(k)$,
and $ \ab{\wE(Z_{k-1},Z_k)}\ab{\vfi_{k-1}''}\leqslant
C\ab{z}r_1(k)$ and Lemma \ref{factorial} allow us to obtain the
expression
\begin{equation}
\vfi_k''=\nu_1(k)z''+\frac{\nu_2(k)}{2}(z^2)''+\frac{(\wE(Z_{k-1},Z_k))''}{\vfi_{k-1}}+C\ab{z}r_1(k).
\label{ce5}\end{equation}
  It follows,  from  (\ref{f4}), that
\begin{equation}
\frac{(\wE(Z_{k-1},Z_k))''}{\vfi_{k-1}}=(\wE(Z_{k-1},Z_k))''
+C\ab{z}r_1(k).\label{ce6}
\end{equation}
Now
$(\wE(Z_{k-1},Z_k))''=\wE(Z_{k-1}'',Z_k)+2\wE(Z_{k-1}',Z_{k}')+\wE(Z_{k-1},Z_k'')$.

\noindent Due to
\[Z_{k-1}'=\ii X_{k-1}\ee^{\ii tX_{k-1}}=z'X_{k-1}(1+\theta
(X_{k-1}-1)\ab{z}/2)=z'X_{k-1}+\theta X_{k-1}(X_{k-1}-1),\]
we obtain
\[
2\wE(Z_{k-1}',Z_k')=2z'\wE(X_{k-1},Z_k')+\theta\wE^{+}_2(X_{k-1},X_k)\ab{z}=
2(z')^2\wE(X_{k-1},X_k)+\theta C\wE^{+}_2(X_{k-1},X_k)\ab{z}.
\]
Similarly, $Z_k=X_kz+\theta X_k(X_k-1)\ab{z}^2/2$
 and
 \[\wE(Z_{k-1}'',Z_k)+\wE(Z_{k-1},Z_k'')=z(\wE(Z_{k-1}'',X_k)+\wE(X_{k-1},Z_k''))+\theta
 C\ab{z}\wE^{+}_2(X_{k-1},X_{k}).\]
Applying (\ref{dui}), we prove
$\wE(Z_{k-1}'',X_k)=z''\wE(X_{k-1},X_k)+\theta
C\wE^{+}_2(X_{k-1},X_k)$.
 Consequently,
 \[(\wE(Z_{k-1},Z_k))''=(z^2)''\wE(X_{k-1},X_k)+\theta
 C\ab{z}\wE^{+}_2(X_{k-1},X_k).\]

 \vspace*{-0.4cm}
\noindent Combining the last estimate with (\ref{ce6}) and (\ref{ce5}), we
 complete the proof of (\ref{vfi2a}). The case $k<6$ is proved exactly by the same
 arguments.
\qed

%%%%%%%%%%%%%%%%%%%%%%%%%%%%%%%%%%%%%%%%%%%%%%%%%%%%%%%%%%%%%%%%%%%%%%%%%%%%%%%%%%%%%%%%%%%%%%%%%%%
%%%DIFFER VFI-G

\newpage
Let $\tilde\vfi_k=\vfi_k\exponent{-\ii t\nu_1(k)},\quad\tilde
g_k=g_k\exponent{-\ii t\nu_1(k)},
\quad\psi =\exponent{-0.1\lambda\sin^2(t/2)}$.

\begin{lemma}\label{figa} Let (\ref{nu12})
hold. Then
\begin{eqnarray*}
\sum_{l=1}^n\ab{\tilde\vfi_l'}&\leqslant& C\lambda\ab{z},\quad
\sum_{l=1}^n\ab{\tilde g_l'}\leqslant C\lambda\ab{z},\quad
\sum_{l=1}^n\ab{\tilde\vfi_l''}\leqslant C\lambda,\\
\sum_{l=1}^n\ab{\tilde g_l''}&\leqslant& C\lambda,\quad
\Ab{\prod_{l=1}^n\tilde\vfi_l-\prod_{l=1}^n\tilde g_l}\leqslant C
R_1\ab{z}^3\psi,\\
\Ab{\Big(\prod_{l=1}^n\tilde\vfi_l-\prod_{l=1}^n\tilde
g_l\Big)'}&\leqslant& C R_1\ab{z}^2\psi,\quad
\Ab{\Big(\prod_{l=1}^n\tilde\vfi_l-\prod_{l=1}^n\tilde
g_l\Big)''}\leqslant C R_1\ab{z}\psi.
\end{eqnarray*}
\end{lemma}

\noindent \Proof The first four estimates follow from Lemmas \ref{fi13} and
\ref{fin} and trivial estimate $\Expect X_{k-1}X_k\leqslant
C_0\nu_1(k)$. Also, using (\ref{nu12}) and (\ref{elk}), we get
\[
\prod_{l=1,l\ne k}^n\exponent{-\lambda_l\sin^2(t/2)}\leqslant C
\prod_{l=1}^n\exponent{-\lambda_l\sin^2(t/2)}\leqslant C\psi^2.
\]
Therefore, by (\ref{elk}) and (\ref{sumi}),
\[
\Ab{\prod_{l=1}^n\tilde\vfi_l-\prod_{l=1}^n\tilde
g_l}=\Ab{\prod_{l=1}^n\vfi_l-\prod_{l=1}^n g_l}\leqslant
\sum_{j=1}^n\ab{\vfi_j-g_j}\prod_{l=1}^{j-1}\ab{g_l}\prod_{l=j+1}^n\ab{\vfi_l}\leqslant
C\psi^2\sum_{j=1}^n\ab{\vfi_j-g_j}\leqslant CR_1\ab{z}^3\psi^2.
\]
From (\ref{nu12}) and trivial estimate $z\ee^{-x}\leqslant 1$, for
$x>0$, we get
\begin{equation*}
\ab{\Gamma_2}\leqslant 0.08\lambda,\quad
\lambda\ab{z}^2\psi\leqslant C. \label{haa1}
\end{equation*}

\vspace*{-0.4cm}
\noindent Therefore,
\begin{eqnarray*}
\Ab{\Big(\prod_{l=1}^n\tilde\vfi_l-\prod_{l=1}^n\tilde
g_l\Big)'}&\leqslant&\sum_{l=1}^n\ab{\tilde\vfi_l'-\tilde
g_l'}\prod_{k\ne l}\ab{\tilde\vfi_k}+\sum_{l=1}^n\ab{\tilde
g_l'}\Ab{\prod_{k\ne l}\tilde\vfi_k-\prod_{k\ne l}\tilde g_k}\\
&\leqslant& C\psi^2[R_1\ab{z}^2+\lambda\ab{z}R_1\ab{z}^3]\leqslant
C\psi R_1\ab{z}^2.
\end{eqnarray*}
The proof of last estimate is very similar and therefore omitted.
\qed

%%%%%%%%%%%%%%%%%%%%%%%%%%%%%%%%%%%%%%%%%%%%%%%%%%%%%%%%%%%%%%%%%%%%%%%%%%%%%%%%%%%%%%%%%%%%%%%%%%%%%%%%%%%%%%%
%%%PROOOFS
%%%%%%%%%%%%%%%%%%%%%%%%%%%%%%%%%%%%%%%%%%%%%%%%%%%%%%%%%%%%%%%%%%%%%%%%%%%%%%%%%%%%%%%%%%%%%%%%%%%%%%%%%%%%%%
\section{Proof of Theorems }

\textbf{Proof of Theorem \ref{NON}.} Hereafter, $x\in\NN$, the set of nonnegative integers. The
beginning of the proof is almost identical to the proof of
Tsaregradsky's inequality. Let $M$ be concentrated on integers.
Then summing up the formula of inversion
 \begin{equation}
M\{k\}=\frac{1}{2\pi}\int_{-\pi}^\pi\w M(t)\ee^{-\ii t k}\dd
t\label{inversion}\end{equation},
 we get
\[\sum_{k=m}^xM\{k\}=\frac{1}{2\pi}\int_{-\pi}^\pi\w
M(t)\frac{\ee^{-\ii t(m-1)}-\ee^{-\ii t x}}{z}\,\dd t.
\]
  If $\ab{\w M(t)/z}$ is bounded,
then as $m\to -\infty$ and by Riemann-Lebesgue theorem, we get
\begin{equation} M(x)=-\frac{1}{2\pi}\int_{-\pi}^\pi\frac{\w
M(t)\ee^{-\ii tx}}{z}\,\dd
t=-\frac{1}{2\pi}\int_{-\pi}^\pi\frac{\w M(t)\ee^{-\ii
t/2}\ee^{-\ii tx}}{2\ii\sin(t/2)}\,\dd t. \label{prad}
\end{equation}
The Tsaregradsky's inequality
\begin{equation}
\ab{M(x)}\leqslant\frac{1}{2\pi}\int_{-\pi}^\pi\frac{\ab{\w
M(t)}}{\ab{z}}\,\dd t \label{Care}
\end{equation}
now follows easily.
\noindent Let next $M=F_n-\G$. Then expressing $\w M(t)$ in powers of $z$, we
get $\w M(t)=\sum_{k=2}a_kz^k$, for some coefficients $a_k$ which
depend on factorial moments of $S_n$. Therefore, $\w
M(\pi)/z(\pi)=\w M(-\pi)/z(-\pi)$. Consequently, integrating  (\ref{prad}) by
parts,
we obtain, for $x\ne\lambda$,
\[
M(x)=-\frac{1}{2\pi}\int_{-\pi}^\pi\frac{\w M(t)\ee^{-\ii
t(\lambda+1/2)}}{2\ii\sin(t/2)}\,\ee^{-\ii t(x-\lambda)}\, \dd t=
\frac{1}{2\pi(x-\lambda)^2}\int_{-\pi}^\pi u''(t)\ee^{-\ii
t(x-\lambda)}\dd t,
\]
where
\[u(t)=\ee^{-(\lambda+1/2)\ii t}\frac{\w M(t)}{2\ii\sin(t/2)}=\frac{\prod_{j=1}^n\tilde\vfi_j-\prod_{j=1}^n\tilde
g_j}{z}.\] Thus, for all $x\in\NN$,
\begin{equation}
(x-\lambda)^2M(x)\leqslant
\frac{1}{2\pi}\int_{-\pi}^\pi\ab{u''(t)}\dd t. \label{u2}
\end{equation}
Using Lemma \ref{figa},  equations (\ref{Care}),  (\ref{u2}) and the trivial estimate
\begin{equation}
\int_{-\pi}^\pi \ab{z}^k\psi(t)\dd t\leqslant
\frac{C(k)}{\lambda^{(k+1)/2}}.\label{trivial}\end{equation}
the proof of (\ref{g}) follows.

All other approximations are compared to compound Poisson measure $\G$ and then the
triangle inequality is applied. We begin from the negative
binomial distribution. Due to the assumptions,
\[ \Gamma_2\leqslant \frac{3}{40}\lambda,\quad
\frac{1-\qubar}{\qubar}=\frac{2\Gamma_2}{\lambda}\leqslant 0.15,\]
see \citep{CeVe13}. Therefore, $\w{NB}(t)\exponent{-\lambda\ii
t}=\exponent{A}$, where
\[A=\lambda{z-\ii
t}+\Gamma_2z^2+\sum_{j=3}^\infty\frac{r}{j}\bigg(\frac{1-\qubar}{\qubar}\bigg)^jz^j=\lambda(z-\ii
t)+\Gamma_2z^2+\theta C\Gamma_2^2\lambda^{-1}\ab{z}^3.\] Moreover,
\[\ab{A'}\leqslant C\lambda\ab{z},\quad\ab{A''}\leqslant
C\lambda,\quad \ab{\ee^A}\leqslant\psi^2.\] Let $B=\lambda(z-\ii
t)+\Gamma_2z^2$ so that $\w\G(t)\exponent{-\lambda\ii
t}=\exponent{B}$ and $u_1(t)= ({\ee^A-\ee^B})/{z} $.
Then
\begin{equation}
\ab{u_1}\leqslant\frac{\ab{\ee^A-\ee^B}}{\ab{z}}\leqslant\psi^2\frac{\ab{A-B}}{\ab{z}}\leqslant
C\psi^2\frac{\Gamma_2^2\ab{z}^2}{\lambda},\quad\int_{-\pi}^\pi\ab{u_1}\dd
t\leqslant C\frac{\Gamma_2^2}{\lambda^2\sqrt{\lambda}}.
\label{carnb}
\end{equation}
Also,
\begin{eqnarray*}
\ab{(\ee^A-\ee^B)''}&\leqslant&\ab{A''}\ab{\ee^A-\ee^B}+\ab{(A')^2}\ab{\ee^A-\ee^B}+
\ab{A''-B''}\ab{\ee^B}+\ab{(A')^2-(B')^2}\ab{\ee^B}\\
&\leqslant& C\psi^2\bigg\{
\lambda\frac{\Gamma_2^2}{\lambda}\ab{z}^3+\lambda^2\ab{z}^2\frac{\Gamma_2^2}{\lambda}\ab{z}^3+
\frac{\Gamma_2^2}{\lambda}\ab{z}+\lambda\ab{z}\frac{\Gamma_2^2}{\lambda}\ab{z}^2
\bigg\}\leqslant C\psi\ab{z}\frac{\Gamma_2^2}{\lambda}.
\end{eqnarray*}
Similarly,
\[
\ab{(\ee^A-\ee^B)'}\leqslant\ab{A'}\ab{\ee^A-\ee^B}+\ab{\ee^B}\ab{A'-B'}\leqslant
C\psi\ab{z}^2\frac{\Gamma_2^2}{\lambda}
\]
and we obtain finally
\begin{equation}
\ab{u_1''}\leqslant
C\psi\frac{\Gamma_2^2}{\lambda},\quad\int_{-\pi}^\pi\ab{u_1''}\dd
t\leqslant C\frac{\Gamma_2^2}{\lambda\sqrt{\lambda}}. \label{unb2}
\end{equation}
Estimates in (\ref{carnb}) and (\ref{unb2}) allow us to write
\[
\bigg(1+\frac{(x-\lambda)^2}{\lambda}\bigg)\ab{\G(x)-\NB(x)}\leqslant
C\frac{\Gamma_2^2}{\lambda^2\sqrt{\lambda}},
\]
which combined with (\ref{g}) proves (\ref{nb}).

\noindent For the proof of
translated Poisson approximation, let $B$ be defined as in above,
\[ T=\lambda (z-\ii t)+(2\Gamma_2+\tilde\delta)(z-\ii t),\quad
D=\lambda (z-\ii t)+(\Gamma_2+\tilde\delta/2)z^2,\] and
\[u_2=(\ee^D-\ee^T)/z,\quad u_3=(\ee^B-\ee^D)/z.\]
Note that, for $\ab{t}\leqslant\pi$, we have
$\ab{t}/\pi\leqslant\ab{\sin(t/2)}\leqslant\ab{t}/2$. Therefore,
arguing similarly as in above, we obtain
\begin{equation}
\int_{-\pi}^\pi\ab{u_2}\dd
t\leqslant\frac{C(\ab{\Gamma_2}+\tilde\delta)}{\lambda\sqrt{\lambda}},\quad
\int_{-\pi}^\pi\ab{u_2''}\dd
t\leqslant\frac{C(\ab{\Gamma_2}+\tilde\delta)}{\sqrt{\lambda}}.
\label{tp1}
\end{equation}
Observe next that
\[u_3=\frac{\ee^B}{z}(\ee^{\tilde\delta z^2/z}-1)=\frac{\ee^B}{z}\int_0^1(\tilde\delta
z^2/2)\ee^{\tau\tilde\delta z^2/2}\dd
\tau=\int_0^1\frac{\tilde\delta z}{2}\ee^{B+\tau\tilde\delta
z^2/2}\dd \tau.
\]
Consequently,
\begin{equation}
\int_{-\pi}^\pi\ab{u_3}\dd t\leqslant C\int_{-\pi}^\pi
\psi^2\tilde\delta\ab{z}\dd t\leqslant
\frac{C\tilde\delta}{\lambda}. \label{tp2}
\end{equation}
Similarly,
\[
u_3''=\frac{\tilde\delta}{2}\int_0^1\ee^{B+\tau\tilde\delta
t}[z''+2z'(B'+\tau\tilde\delta
zz')+z(B''+\tau\tilde\delta(zz')')+z(B'+\tau\tilde\delta
zz')^2]\dd\tau
\]
and using $\tilde\delta\leqslant 1\leqslant\lambda$, we get
\[\ab{u_3''}\leqslant
C\psi^2\tilde\delta(1+\lambda\ab{z}+\tilde\delta\ab{z}+\ab{z}(\lambda\ab{z}+\tilde\delta\ab{z})^2)
\leqslant C\tilde\delta\psi\sqrt{\lambda}.
\]
Consequently,
\[\int_{-\pi}^\pi\ab{u_3''}\dd t\leqslant C\tilde\delta.\]
Combining the last estimate, the inequalities in (\ref{tp1}), (\ref{tp2}) and the
estimate for $\w G=\ee^B$, the result in (\ref{tp}) is proved.

\noindent For binomial approximation, note first that
\vspace*{-0.4cm}
\begin{eqnarray*}\ee^{-\lambda\ii t}\w \BI&=&\ee^E,\quad E=\lambda(z-\ii
t)+\Gamma_2z^2+z^2\theta\frac{50\Gamma_2^2}{21\lambda^2}\varepsilon+\theta\frac{5N\pbar^3\ab{z}^3}{9},\\
\pbar&\leqslant&\frac{50\ab{\Gamma_2}}{21\lambda}<\frac{1}{5},\quad\ab{\Gamma_2}\leqslant
0.08\lambda,\quad\ab{N\pbar^3}\leqslant
C\frac{\Gamma_2^2}{\lambda},
\end{eqnarray*}
see \citep{CeVe13}. Let

\vspace*{-0.4cm}
\[
L=\lambda(z-\ii
t)+\Gamma_2z^2+z^2\theta\frac{50\Gamma_2^2}{21\lambda^2}\epsilon,
\quad u_4=(\ee^L-\ee^E)/z,\quad u_5=(\ee^B-\ee^L)/z.
\]
Next,
\[u_5=\int_0^1\ee^B
z\Exponent{\tau
z^2\theta\frac{50\Gamma_2^2}{21\lambda^2}\epsilon}\theta\frac{50\Gamma_2^2}{21\lambda^2}\epsilon\dd\tau.
\]
Now the proof is practically identical to that  of (\ref{tp})
and is, therefore, omitted.

\noindent The proofs of (\ref{pi}) and (\ref{pi2}) are also very similar
 and use the facts
\begin{eqnarray*}
\frac{\ee^B-\ee^{-\lambda\ii
t}(\w\Pi+\w\Pi_1)}{z}&=&\int_0^1(1-\tau)\Gamma_2^2z^3\exponent{\lambda(z-\ii
t)+\tau\Gamma_2z^2}\dd\tau,\\
 \frac{\ee^B-\ee^{-\lambda\ii
t}\w\Pi}{z}&=&\int_0^1\Gamma_2z\exponent{\lambda(z-\ii
t)+\tau\Gamma_2z^2}\dd\tau.  \rlap{$\qquad \Box$}
\end{eqnarray*}

\noindent \textbf{Proof of Theorem \ref{NONL}.}
 Let $M$ be a measure concentrated on integers and $\w
 M(t)=\sum_{k=1}^\infty M\{k\}\ee^{\ii tk}$.
 Then from formula  (\ref{inversion}) of
 inversion,  we get
\[
\ab{M\{x\}}\frac{1}{2\pi}\leqslant\int_{-\pi}^\pi\ab{\w M(t)}\dd
t.
\]

\vspace*{-0.4cm}
Moreover, integrating (\ref{inversion}) by parts, we obtain
\[
(x-\lambda)^2\ab{M\{x\}}\leqslant\frac{1}{2\pi}\int_{-\pi}^\pi\ab{(\w
M(t)\exponent{-\lambda\ii t})''}\dd t.
\]
The rest of the proof is a simplified version of the proof of Theorem
\ref{NON} and hence omitted. \qed

\vspace*{0.4cm}
\noindent {\bf Acknowledgments}.  The  authors are grateful to Dr. Sriram for inviting us to contribute this
article for Dr. Koul's Festschrift and to the referee for several helpful comments.

%%%%%%%%%%%%%%%%%%%%%%%%%%%%%%%%%%%%%%%%%%%%%%%%%%%%%%%%%%%%%%%%%%%%%%%%%%%%%%%%%%%%%%%%%%%%%%%%%%%%%%%%%%%%%%%%%%%%%%%%%%%%%%

%%%%%%%%%%%%%%%%%%%%%%%%%%%%%%%%%%%%%%%%%%%%%%%%%%%%%%%%%%%%%%%%%%%%%%%%%%%%%%%%%%%%%%%%%%%%%%%%%%%%%%%%%%%%%%%%%%%%%%%%
%%%%% REFERENCES
%%%%%%%%%%%%%%%%%%%%%%%%%%%%%%%%%%%%%%%%%%%%%%%%%%%%%%%%%%%%%%%%%%%%%%%%%%%%%%%%%%%%%%%%%%%%%%%%%%%%%%%%%%%%%%%%%%%%%%%
%%%%%%%%%%%%%%%%%%%%%%%%%%%%%%%%%%%%%%%%%%%%%%%%%%%%%%%%%%%%%%%%%%%%%%%%%%%%%%%%%%%%%%%%%%%%%%%%%%%%%%%%%%%%%%%%%%%%%%%%

\end{document}